\documentclass[12pt]{article}
\usepackage[T1]{fontenc}
\usepackage[english]{babel}
\usepackage{amssymb,url,xspace,amsmath,amsthm,eucal,mathrsfs,eufrak,array,epsfig,psfrag,graphicx}
\usepackage[all]{xy}

\newtheorem{thm}{Theorem}[section]
\newtheorem{lem}[thm]{Lemma}
\newtheorem{prop}[thm]{Proposition}

\theoremstyle{definition}

\newtheorem{defn}[thm]{Definition}
\newtheorem{ex}[thm]{Example}
\def \srec  {\mathop{\rm srec}\limits}
\def \rec  {\mathop{\rm rec}\limits}

\newenvironment{preuve}[1][Proof.]{\begin{trivlist}
\item[\hskip \labelsep {\bfseries #1}]}{\fin \end{trivlist}}

\def \ee  {\mathop{\rm sup}\limits}
\def \N {\mathbb{N}}
\def \S {\mathcal{S}}

\def \C {\mathcal{C}}

\def \fin {\hfill $\blacksquare$}

\def \N {{\mathbb N}}

\def \S {{\mathfrak S}}

\def \Card {\mathop{{\rm Card } }\nolimits}

\begin{document}

\vskip1cm

\centerline{\LARGE Asymptotic behavior of permutation records}

\vspace{8mm}

\centerline{\large Igor Kortchemski\footnote{E-mail address:
\texttt{igor.kortchemski@ens.fr}}}

\vspace{8mm}

\centerline{\large \it \'{E}cole Normale Supérieure, 75005 Paris,
France}

 \vspace{3mm}
\vspace{1cm}

\begin{abstract}
We study the asymptotic behavior of two statistics defined on the
symmetric group $\S_n$ when $n$ tends to infinity: the number of
elements of $\S_n$ having $k$ records, and the number of elements of
$\S_n$ for which the sum of the positions of their records is $k$.
We use a probabilistic argument to show that the scaled asymptotic
behavior of these statistics can be described by remarkably simple
functions.
\end{abstract}

\section{Introduction}

In this paper we will study records, also known in the literature as
\emph{left to right maxima}, \emph{outstanding elements} or
\emph{strong records}. By definition, a \emph{record} of a
permutation $\sigma = a_{1} \ldots a_{n} \in \S_n$ is a number
$a_{j}$ such that $a_{i} < a_{j}$ for all $i < j$. The study of
records has been initiated by Rényi who proved that the number of
elements of $\S_n$ with exactly $m$ cycles in their cycle
decomposition is equal to the number of elements of $\S_n$ having
exactly $m$ records, the latter being given by $c(n,m)$, the
unsigned Stirling number of the first kind \cite{Renyi}. The
asymptotic behavior of these numbers has first been studied by
Jordan \cite{Jordan} who showed that for fixed $m$ and large $n$,
$c(n,m)/(n-1)!\sim(\ln n + \gamma)^{m-1}/(m-1)!$, with $\gamma$
being the Euler constant. Moser and Wyman \cite{MoserWyman}
considered three overlapping regions of the $(n,m)-$plane (for
$n\geq m$) and obtained asymptotic formulae in each case. Wilf
\cite{Wilf} gave an explicit asymptotic expansion of $c(n,m)$ with
$m$ fixed in terms of powers of $n$ and $\ln n$. Finally, Temme
\cite{Temme} obtained an asymptotic formula for $c(n,m)$ in the
limit $n\rightarrow \infty$, uniformly for $1\leq m\leq n$
(including the case $m\sim n$). Hwang \cite{Hwang} established an
explicit asymptotic expansion of $c(n,m)$ valid for $m
=\mathcal{O}(\ln n)$. Later, Chelluri, Richmond and Temme
\cite{Temme2} introduced the generalized Stirling numbers of the
first kind and studied their asymptotic behavior. In particular,
they re-derived Temme's previous result and showed that their
results agree with those of Moser and Wyman.

 Asymptotic properties of record statistics have been studied by Wilf who obtained
the following results \cite{Wilf1995}:

(i) for fixed $r$ the average value of the $r$th records, over all
permutations of $\S_n$ that have that many, is asymptotic to $(\ln
n)^{r-1}/(r - 1)!$ when $n\rightarrow+\infty$,

(ii) the average value of a permutation $\sigma$ at its $r$th
record, among all permutations that have that many, is asymptotic to
$(1-1/2^r)n$ when $n\rightarrow+\infty$,

(iii) Let $1 < j_1 < j_2 < \cdots < j_m$ be fixed integers, and
suppose that we attach a symbol $s(j) =$ `Y' or `N' to each of these
$j$'s. Then the probability $P$ that a permutation of $\S_n$ does
have a record at each of the $j_v$ that is marked `Y', and does not
have a record at any of those that are marked `N',
is:\begin{equation}\label{eqn:Wilf}
P=\prod_{s(j_v)=\textrm{`N'}}\left(1-\frac{1}{j_v}\right)\prod_{s(j_v)=\textrm{`Y'}}\frac{1}{j_v}\,.\end{equation}

Myers and Wilf generalized the notion of records \cite{MW}. They
studied strong and weak records defined on multiset permutations and
words (a \emph{weak record} of a word $w_1,\ldots,w_n$ is a term
$w_{j}$ such that $w_{i} \leq w_{j}$ for all $i < j$ \cite{Geo}).
For multiset permutations, they derived the generating function for
the number of permutations of a fixed multiset $M$ which contain
exactly $r$ strong (respectively weak) records. They also obtained
the generating function of the probability that a randomly selected
permutation of $M$ has exactly $r$ strong records. This gives the
average number of strong (respectively weak) records among all
permutations of $M$.

The notion of records has been extended to random variables (for a
survey of some results see \cite{Glick}). The asymptotic behavior of
the two statistics, `position' and `value' of the $r$th record, have
been studied in \cite{Geo0,Geo,Geo2} for a sequence of i.i.d random
variables which follow the geometric law of parameter $p$, which
approach the model of random permutations in the limit $p\rightarrow
0$. For other interesting results concerning record statistics see
\cite{MW}.

\bigskip

In this article, we re-derive Wilf's result (\ref{eqn:Wilf}) by
using of a \emph{probabilistic argument}. This will allow us to
perform a study of the asymptotic behavior of another statistic of
records: the number of permutations of length $n$ having $k$ records
in the limit $n\rightarrow+\infty$ with the ratio $k/n$ fixed. In
this limit the \emph{re-scaled} number of records $k/n$ takes values
on the interval $[0,1]$. We also introduce the new statistic for a
permutation called `sum of the positions of its records'. We find
the asymptotic behavior of the number of permutations of length $n$
for which the sum of the positions of their records is $k$ in the
limit $n\rightarrow+\infty$ with the ratio $k/(\frac{n(n+1)}{2})$
fixed. More precisely, we show the following results:

\bigskip

(I) Let $c(n,k)$ be the number of permutations of length $n$ having
$k$ records. Its generating function is given by
$q(q+1)\cdots(q+n-1)$, so that $c(n,k)$ is the coefficient of $q^k$
in the power expansion. For $n \geq 1$ and $x \in [0,1]$ define the
function $f_n$ by:

\begin{equation}\label{defn:fn}f_n(x) = \begin{cases}
       c\left(n,\left[nx\right]\right)  & \text{if } x \geq \frac{1}{n}, \\
       c(n,1)  & \text{otherwise,}
       \end{cases}\end{equation}
where $[x]$ stands for the integer part of $x$. Then, when $n$ tends
to infinity, the sequence of functions
$\{\frac{\ln(f_n)}{n\ln(n)};n\in \N^+\}$ converges uniformly with
respect to $x$ on the interval $[0,1]$ to the function $x\mapsto1-x$
with an accuracy $\mathcal{O}\left(\frac{1}{\ln n}\right)$. In other
words, there exists a constant $C$ such that for all integer $n\geq
2$:
$$\ee_{x\in[0,1]}\left|\frac{\ln(f_n(x))}{n\ln n}-(1-x)\right|\leq \frac{C}{\ln n}\,.$$

\bigskip

(II) Let $\C(n,k)$ be the number of permutations of length $n$ for
which the sum of the positions of their records is $k$. Its
generating function is given by $q(q^2+1)(q^3+2)\cdots(q^n+n-1)$, so
that $\C(n,k)$ is the coefficient of $q^k$ in the power expansion.
For $n\geq 1$ and $x \in [0,1]$ define the function $\phi_n$
by\footnote{This definition is motivated by the fact that
$\C(n,k)=0$ if and only if $k=2$ or $k=\frac{n(n+1)}{2}-1$}:

\begin{equation}\label{defn:phin}\phi_n(x) = \begin{cases}
        \C(n,1)=(n-1)!  & \text{if } x<\frac{6}{n(n+1)}\\
        \C(n,\frac{n(n+1)}{2})=1 & \text{if }
        x\geq 1-\frac{2}{n(n+1)}\\
       \C\left(n,\left[\frac{n(n+1)}{2}x\right]\right) & \text{otherwise.}
       \end{cases}\end{equation}
Then, when $n$ tends to infinity, the sequence of functions
$\{\frac{\ln(\phi_n)}{n\ln(n)};n\in \N^+\}$ converges uniformly with
respect to $x$ on the interval $[0,1]$ to the function
$x\mapsto\sqrt{1-x}$ with an accuracy $\mathcal{O}\left(\frac{1}{\ln
n}\right)$. In other words, there exists a constant $\widetilde{C}$
such that for all integer $n$:
$$\ee_{x\in[0,1]}\left|\frac{\ln(\phi_n(x))}{n\ln n}-\sqrt{1-x}\right|\leq \frac{\widetilde{C}}{\ln n}.$$

\bigskip
It is important to note that the functions $f_n$ and $\phi_n$ are
defined on the same interval $[0,1]$.

\bigskip
The paper is organized as follows. In section $2$ we derive some
results concerning records which will be used in the rest of the
paper. In section $3$ we prove assertion (I) and in section $4$ we
show assertion (II). In the appendix we prove that (I) is consistent
with Temme's result \cite{Temme} previously mentioned.

\section{Notations and preliminary results}

We endow the symmetric group $\S_n$ on a set of $n$ elements with
the uniform law. We begin with some useful definitions.

\begin{defn}Let $\sigma= a_{1} \ldots a_{n} \in \S_n$. Recall that a \emph{record} of $\sigma$
 is a number $a_{j}$ such that $a_{i} < a_{j}$ for all $i < j$. We define $\textrm{rec}(\sigma)$ as the
 number of records of $\sigma$. The generating function of this statistic is:
$$T_n(q)=\sum_{\sigma \in \S_n}q^{\rec(\sigma)}.$$
Likewise we define $\textrm{srec}(\sigma)$ as the sum of the
positions of all records of $\sigma$. The generating function of
this statistic:
$$P_n(q)=\sum_{\sigma \in \S_n}q^{\srec(\sigma)}.$$
Let $X_k(\sigma)$ be the random variable which equals $1$ if $k$ is
a position of a record of $\sigma$ and  $0$ otherwise.\end{defn}

\begin{ex}The permutation of length $8$, $\sigma=4,7,5,1,6,8,2,3$ (i.e.
$\sigma$ sends $1$ to $4$, $2$ to $7$ etc.) has $3$ records: $4$,$7$
and $8$ so that $\textrm{rec}(\sigma)=3$ and
$\textrm{srec}(\sigma)=1+2+6=9$.\end{ex}

Our work relies on the following proposition, first proved by Rényi
\cite{Renyi}.

\begin{prop}\label{prop:indep}The random variables $X_1, X_2,\ldots,X_n$
are independent. Moreover,
$\mathbb{P}(X_k=1)=\frac{1}{k}$.\end{prop}

\begin{preuve}This proposition is a consequence of proposition 1.3.9
and corollary 1.3.10 of \cite{PR,ENUM}. It comes from the following
remark. For a permutation $\sigma=a_1\ldots a_n$ and an integer $n$
such that $1\leq i\leq n$, define:
$$ r_i(\sigma) =\Card\{j\colon j<i,\ a_j>a_i\}.$$
Then the mapping which sends a permutation $\sigma=a_1\ldots a_n$ on
the $n$-tuple $(r_1(\sigma),\ldots,r_n(\sigma))$ is as a bijection
between $\S_n$ and the $n$-tuples $(r_1,\ldots,r_n)$ such that
$0\leq r_i\leq i-1$ for all $i$. Furthermore, $r_i(\sigma)=0$ if,
and only if, $a_j$ is a record of $\sigma$.
\end{preuve}

\begin{ex}By this bijection, the image of $4,7,5,1,6,8,2,3$ is
$\{0,0,1,3,1,0,5,5\}$.\end{ex}

\begin{defn}\label{defCnk}Define
$\mathcal{C}(n,k)$ as the number of elements of $\S_n$ for which the
sum of the positions of their records is $k$.
\end{defn}

Using proposition \ref{prop:indep} one can show the following
results.

\begin{prop}The generating functions of the statistics `rec' and `srec' are:
\begin{align}\label{rec}T_n(q)&=\sum_{\sigma \in
\S_n}q^{\rec(\sigma)}=q(q+1)\cdots(q+n-1)&&=\sum_{k=0}^n c(n,k)q^k\\
\label{srec} P_n(q)&=\sum_{\sigma \in
\S_n}q^{\srec(\sigma)}=q(q^2+1)(q^3+2)\cdots(q^n+n-1)&&=\sum_{k=1}^{\frac{n(n+1)}{2}}\mathcal{C}(n,k)q^k.
\end{align}
\end{prop}

The generating function (\ref{rec}) is well known \cite{ENUM}.

\section{Asymptotic behavior of the coefficients $c(n,k)$}

Let us first examine the asymptotic behavior of $c(n,k)$ when $n$
tends to infinity and $k/n$ is fixed. In this limit, it is not
obvious that the coefficients $c(n,k)$ have a well defined
asymptotic behavior. It is convenient to introduce the new scaling
variable $x =k/n$ which takes values in the interval $[0,1]$. We
introduce a new function $f_n$ as:
\begin{equation}\label{ee1}
f_n(x)=c\left(n,nx\right)\,.
\end{equation}

 Note that when $x$ is of the form
$k/n$ the two relations (\ref{defn:fn}) and (\ref{ee1}) coincide.
For a fixed integer $n$, the variable $x$ takes rational values.
Moreover in the limit $n \rightarrow \infty$ they run through a
dense subset of $[0,1]$. One may wonder whether the function
$f_n(x)$ is well defined in this limit. And if so, what is the limit
function?

For this purpose we extend the function $f_n$ to the whole interval
$[0,1]$ as in Eq.(\ref{defn:fn}). We now state the theorem which
describes the asymptotic behavior of $c(n,k)$.

\begin{thm}\label{thm1}The sequence of
functions $\{\frac{\ln(f_n)}{n\ln(n)};n\in \N^+\}$ converges
uniformly when $n$ tends to infinity with respect to $x$ on the
interval $[0,1]$ to the function $x\mapsto1-x$ with an accuracy
$\mathcal{O}\left(\frac{1}{\ln n}\right)$. In other words, there
exists a constant $C$ such that for all integer $n$:
\begin{equation}\label{eqn:th1}\ee_{x\in[0,1]}\left|\frac{\ln(f_n(x))}{n\ln(n)}-(1-x)\right|\leq \frac{C}{\ln
n}.\end{equation}\end{thm}

The proof relies on a combinatorial and probabilistic interpretation
of the coefficients $c(n,k)$. In the rest of this section we
consider $n$ to be an integer.

\begin{lem}\label{lemme:formule}Let $\mathbb{P}(\rec=k)=\frac{c(n,k)}{n!}$. Then:
\begin{equation}\label{formule}
\mathbb{P}(rec=k)=\sum_{\substack{v_1<v_2<\cdots<v_k\leq n\\
v_1=1}}\frac{1}{v_1v_2\cdots
v_k}\left(1-\frac{1}{v_{k+1}}\right)\cdots\left(1-\frac{1}{v_{n}}\right)
\end{equation}
under the additional conditions $v_{k+1}<\cdots<v_n\leq n$ and $v_i
\neq v_j$ for $i \neq j$.
\end{lem}

\begin{preuve}Choosing a permutation with $k$ records is equivalent to
choosing the positions of its $k$ records as
$1=v_1<v_2<\cdots<v_k\leq n$. By proposition \ref{prop:indep}, a
position $v_i$ is chosen as a record position with probability
$\frac{1}{v_i}$ and is not chosen as a record position with
probability $1-\frac{1}{v_i}$. Furthermore this choice is
independent for each position. This yields the result of the lemma.
\end{preuve}
This gives a probabilistic proof of Wilf's result (\ref{eqn:Wilf})
mentioned in the Introduction.

\begin{lem}\label{lem:encadrement}Let $x\in[\frac{1}{n},1]$ and $k=[nx]$. The following double
inequality holds:
$$\frac{(n-[nx])!}{n(n!)}\leq \mathbb{P}(rec=k) \leq 2^n \frac{1}{[nx]!}$$
\end{lem}

\begin{preuve}This is a consequence of lemma \ref{lemme:formule} after taking into account the properties:
\begin{enumerate}
\item[(i)] the number of $k$-tuples $(v_1,\ldots,v_k)$ such that
$v_1<v_2<\cdots<v_k \leq n$, $v_1=1$, $v_i \neq v_j$ for $i \neq j$
and $v_{k+1}<\cdots<v_n\leq n$ is less than $2^n$,
\item[(ii)]we have:
$$\frac{1}{n}=\left(1-\frac{1}{2}\right)\cdots\left(1-\frac{1}{n}\right)\leq\left(1-\frac{1}{v_{k+1}}\right)\cdots\left(1-\frac{1}{v_{n}}\right)\leq 1,$$
\item[(iii)] and for integers $1=v_1<v_2<\cdots<v_k\leq n$:
$$k! \leq v_1v_2\cdots v_k \leq \frac{n!}{(n-k)!}.$$
\end{enumerate}
\end{preuve}

All the essential ingredients have now been gathered, and we turn to
the proof of theorem \ref{thm1}.

\medskip

\textbf{Proof of theorem \ref{thm1}}. Using the same variables as in
lemma \ref{lem:encadrement} and the function $f_n$ defined in
(\ref{defn:fn}), one obtains:
\begin{equation}\label{ineg1}\frac{\ln\left((n-[nx])!\right)}{n \ln n}-\frac{1}{n}\leq
\frac{\ln(f_n(x))}{n\ln(n)}\leq \frac{\ln 2}{\ln n}+\frac{\ln n!}{n
\ln n}-\frac{\ln([nx]!)}{n\ln n}.\end{equation}

Define $M_n(x)=\ln n\left|\frac{\ln(f_n(x))}{n\ln(n)}-(1-x)\right|$.
Then:
$$\ee_{x\in[0,1]}M_n(x)\leq\ee_{x\in[0,\frac{1}{n}]}M_n(x)+\ee_{x\in[\frac{1}{n},1]}M_n(x)$$

Using $c(n,1)=(n-1)!$ (see Eq.(\ref{rec})) and lemma 3.3 one
obtains:

$$\ee_{x\in[0,1]}M_n(x)\leq\ee_{x\in[0,\frac{1}{n}]}\ln n\left|\frac{\ln(n-1)!}{n \ln n}-(1-x)\right|+\ee_{x\in[\frac{1}{n},1]}\ln n\left|\frac{\ln(f_n(x))}{n\ln(n)}-(1-x)\right|$$

 Stirling's formula ($\ln n! = n \ln n + \mathcal{O}(n)$) shows that the first term is bounded when $n$ tends to infinity.
 Now denote $B_n$ the second term of the right-hand side of the above inequality. By Eq.(\ref{ineg1}) there exist constants $C_1$ and $C_2$ such that:

$$B_n\leq C_1+\ee_{x\in[\frac{1}{n},1]}\left|\frac{\ln n!}{n}-\frac{\ln([nx]!)}{n}-(1-x)\ln
n\right|+\ee_{x\in[\frac{1}{n},1]}\left|\frac{\ln\left((n-[nx])!\right)}{n}-(1-x)\ln
n\right|\leq C_2,$$ where in the second inequality we used
Stirling's formula. This concludes the proof.\fin

\section{Asymptotics of the coefficients $\mathcal{C}(n,k)$}

We now investigate the asymptotic behavior of the coefficients
$\mathcal{C}(n,k)$ (see definition \ref{defCnk}) for large $n$. When
$n$ is a fixed integer, the coefficients $\mathcal{C}(n,k)$ (with $1
\leq k \leq \frac{n(n+1)}{2}$) are positive integers. As before, we
are interested in their asymptotic behavior in the limit
$n\rightarrow+\infty$ with the ratio $k/(\frac{n(n+1)}{2})$ fixed.
To this end we introduce the scaling variable $x =2k/(n(n+1))$ which
takes values in the interval $[0,1]$ as well as the function
$\phi_n$:
\begin{equation}\label{e1.1}
\phi_n(x)=\C\left(n,\frac{n(n+1)}{2}x\right)\,.
\end{equation} Then, the problem reduces to finding the asymptotic behavior of $\phi_n(x)$ in the limit
$n\rightarrow \infty$. For this purpose we extend the function
$\phi_n$ to the whole interval $[0,1]$ as in Eq.(\ref{defn:phin}).
We now state the main theorem which describes the asymptotic
behavior of the coefficients $\mathcal{C}(n,k)$.

\begin{thm}\label{thm2}
 The sequence of functions
$\{\frac{\ln(\phi_n)}{n\ln(n)};n\in \N^+\}$ converges uniformly when
$n$ tends to infinity with respect to $x$ on the interval $[0,1]$ to
the function $x\mapsto\sqrt{1-x}$ with an accuracy
$\mathcal{O}\left(\frac{1}{\ln n}\right)$. In other words, there
exists a constant $\widetilde{C}$ such that for all integer $n$:
$$\ee_{x\in[0,1]}\left|\frac{\ln(\phi_n(x))}{n\ln(n)}-\sqrt{1-x}\right|\leq \frac{\widetilde{C}}{\ln n}.$$\end{thm}

To illustrate this theorem we plot the two functions $x \mapsto
\sqrt{1-x}$ and $\psi_n=\frac{\ln(\phi_n)}{n\ln(n)}$ for different
values of $n$ (for $n=50$ in fig. 1 and for $n=150$ in fig. 2)
obtained with \verb"Mathematica". In agreement with theorem
\ref{thm2}, $\psi_n$ converges to the function $\sqrt{1-x}$. We also
plot the function:
$$ \tau: \quad n \longmapsto \ln n \cdot
\ee_{x\in[0,1]}\left|\frac{\ln(\phi_n(x))}{n\ln(n)}-\sqrt{1-x}\right|$$
in fig. 3 for $n=2,3,\ldots,50$. In agreement with theorem
\ref{thm2}, this function is bounded by a constant $\widetilde{C}$.

\begin{figure}
\begin{minipage}[t]{.4\linewidth}
    \begin{center}
       \includegraphics[scale=0.7]{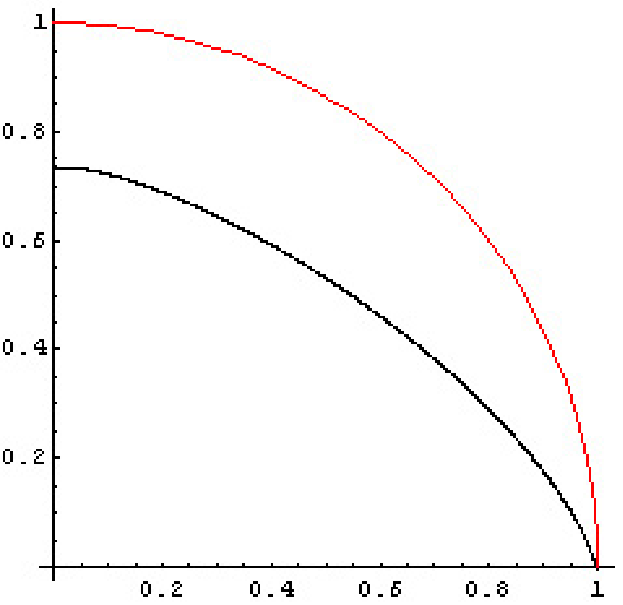}
       \caption{Graph of $\psi_{50}$ and $\sqrt{1-x}$}
       \label{nbabo}
    \end{center}
\end{minipage}
\hfill
\begin{minipage}[t]{.4\linewidth}
    \begin{center}
       \includegraphics[scale=0.7]{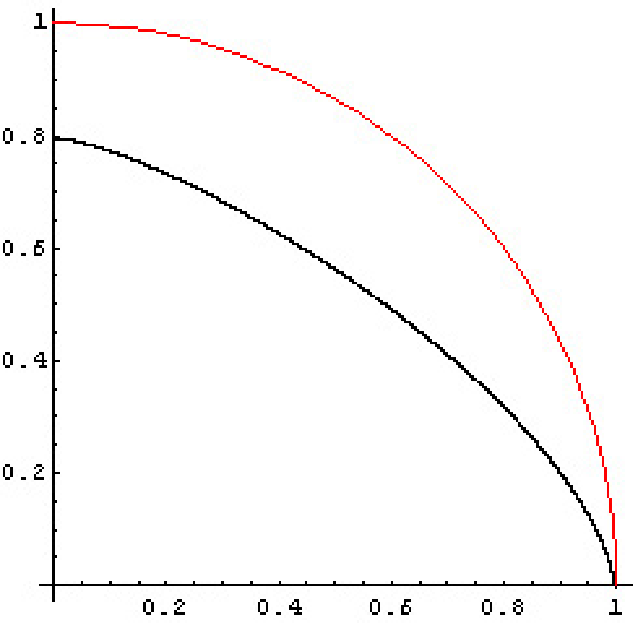}
       \caption{Graph of $\psi_{150}$ and $\sqrt{1-x}$}
       \label{croissnbabo}
    \end{center}
\end{minipage}
\end{figure}

\begin{figure}
    \begin{center}
       \includegraphics[scale=0.7]{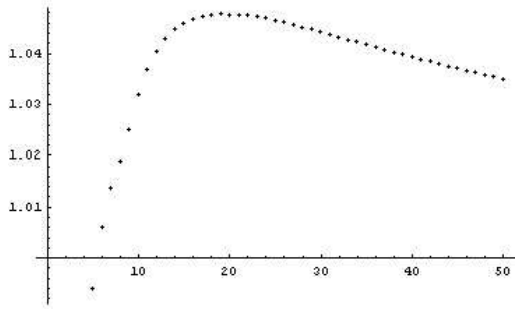}
       \caption{Graph of the function $\tau$ for $n=2,\ldots,50$}
    \end{center}
\end{figure}

\bigskip
The proof of this theorem is based on a combinatorial and
probabilistic interpretation of the coefficients $\C(n,k)$ which
gives us lemma \ref{lemme:formule2}.

\begin{lem}\label{lemme:formule2}Let $\mathbb{P}(\srec=k)=\frac{\mathcal{C}(n,k)}{n!}$. Then:
\begin{equation}\label{formule1}
\mathbb{P}(srec=k)=\sum_{\substack{v_1+v_2+\cdots+v_r=k\\
v_1=1,r\leq n}}\frac{1}{v_1v_2\cdots
v_r}\left(1-\frac{1}{v_{r+1}}\right)\cdots\left(1-\frac{1}{v_{n}}\right)
\end{equation}
under the additional conditions $v_1<v_2<\cdots<v_r\leq n$,
$v_{r+1}<\cdots<v_n \leq n$ and $v_i \neq v_j$ for $i \neq j$.
\end{lem}
\begin{preuve}The proof is similar to the one of lemma
\ref{lemme:formule}. Choosing an element of $\S_n$ for which the sum
of the positions of its records is $k$, is choosing the position of
its records $1=v_1<v_2<\cdots<v_r\leq n$ such that
$v_1+v_2+\cdots+v_r=k$. By proposition \ref{prop:indep}, a position
$v_i$ is chosen as a record position with probability
$\frac{1}{v_i}$ and is not chosen as a record position with
probability $1-\frac{1}{v_i}$. Furthermore this choice is
independent for each position. This yields the result of the lemma.
\end{preuve}

\begin{defn}
We say that the $r$-tuple $(v_1,\ldots,v_r)$ satisfies the
conditions $(C_{k,n})$ if:
\begin{equation}\label{Ckn}(C_{k,n}): \quad \quad v_1=1,\, r\leq n,\,
v_1<v_2<\cdots<v_r\leq n \textrm{ and }
v_1+v_2+\cdots+v_r=k.\end{equation}
\end{defn}

Note that there exists an $r$-tuple $(v_1,\ldots,v_r)$ satisfying
the conditions $(C_{k,n})$ if and only if $\C(n,k)>0$, that is
$k\neq2$ and $k\neq \frac{n(n+1)}{2}-1$. We now isolate the greatest
term in the sum of formula (\ref{formule1}). This motivates the
following definition.

\begin{defn}\label{def:min} Let $k,n$ be integers such that $1 \leq k \leq \frac{n(n+1)}{2}$ , $k\neq2$ and $k\neq \frac{n(n+1)}{2}-1$. Define:
\begin{equation}\label{defmin}m(n,k)=\min_{\substack{r\leq n,(v_1,\ldots,v_r)\\\mbox{
\tiny{satisfies $(C_{k,n})$}}}} v_1v_2\cdots v_r.\end{equation}
\end{defn}

This minimum gives us an inequality satisfied by the coefficients:

\begin{prop}\label{prop:bounds}We have:
\begin{equation}\frac{1}{n m(n,k)}\leq \mathbb{P}(srec=k)
\leq \frac{2^n}{m(n,k)}.\end{equation}
\end{prop}
\begin{preuve}
Note that the total number of $r$-tuples ($1\leq r\leq n$)
satisfying the conditions $(C_{k,n})$ is less then $2^n$, so that:
$$\frac{1}{n m(n,k)}=\frac{1}{m(n,k)}\left(1-\frac{1}{2}\right)\cdots\left(1-\frac{1}{n}\right)\leq \mathbb{P}(srec=k)
\leq \frac{2^n}{m(n,k)}.$$ This concludes the proof.
\end{preuve}

 We will have to study three cases: $3\leq k
\leq n$, $n \leq k < \frac{n(n-1)}{2}$ and $\frac{n(n-1)}{2}\leq k
\leq \frac{n(n+1)}{2}$. For each case, we will find either the
expression for the minimum $m(n,k)$ defined in (\ref{defmin}) or a
lower and upper bound for $m(n,k)$. These expressions will be useful
in finding the asymptotic behavior. In the following, $n$ will be
considered as an integer greater than $3$ and $k$ as an integer.

\subsection{Case $3 \leq k \leq n$}
This case is the easiest one as shows the following lemma.

\begin{prop}\label{prop1}Let $3 \leq k \leq n$. Then:\begin{equation}\frac{n!}{(k-1) n}\leq
\C(n,k) \leq \frac{2^n n!}{k-1}.\end{equation}
\end{prop}

\begin{preuve}
One can verify that the minimum $m(n,k)$, defined in (\ref{defmin}),
is realized by the $r$-tuple $(v_1,\ldots,v_r)$ for $r=2$, when
$v_1=1$ and $v_2=k-1$. In this case, $m(n,k)=k-1$, and proposition
\ref{prop:bounds} yields the result.
\end{preuve}

It is important to note that when $n+1\leq k$ this argument cannot
be applied anymore. We would like to take $v_2=k-1$, this is
impossible since the conditions (\ref{Ckn}) require $v_2\leq n$.

\subsection{Case $n+1\leq k < \frac{n(n-1)}{2}$}
Let $k$ be an integer such that $n+1\leq k < \frac{n(n-1)}{2}$. In
this case, one expects the minimum $m(n,k)$ to be realized when most
of the records are at the last positions:

\begin{lem}\label{lem0} Let $(v_1,\ldots,v_r)$ be an $r$-tuple which
realizes the minimum $m(n,k)$. Let $i_0=i_0(n,k)$ be the greatest
integer such that: \begin{equation}\label{eqn:i0}k-1\geq
n+(n-1)+\cdots+(n-i_0).\end{equation} Then for $0\leq i \leq
i_0(n,k)$, $n-i$ is equal to one of the $v_j$.
\end{lem}

\begin{preuve}
For the sake of contradiction assume that $i$ is the smallest
integer such that $n-i$ does not appear in $v_1,\ldots,v_r$. Then
either $1<v_2<n-i$, or
$$v_1+v_2+\cdots+v_r=1+(n-i+1)+\cdots+(n-1)+n \leq
1+(k-1)-(n-i)<k,$$ which contradicts the conditions (\ref{Ckn}). So
let $j\geq 2$ be the greatest integer such that $v_j<n-i$. The
desired contradiction will arise if we find an $r$-tuple such that
the product of its elements will be less than $m(n,k)$, therefore
contradicting the minimality of $m(n,k)$. A few cases have to be
studied. It is important to remember that the $v_i$ have to be all
different.

\quad \mbox{\emph{Case 1: $v_j \neq n-i-1$}.} If there exists $l >
1$ such that $1<v_l\leq (n-i)-v_j$, we delete $v_l$ and $v_j$ from
our $r$-tuple and replace them by $v_l+v_j$. This contradicts the
minimality of $m(n,k)$ since $v_lv_j$>$v_l$+$v_j$. Moreover if $j=2$
then:
$$1+v_j+(n-i+1)+\cdots+n<1+(n-i)+(n-i+1)+\cdots+n\leq k,$$
which contradicts (\ref{Ckn}). Hence $1<v_2<v_j$ and
$v_2>(n-i)-v_j\geq 2$, implying $v_2-1\geq 2$. Replace $v_2$ by
$v_2-1$ and $v_j$ by $v_j+1$. This contradicts the minimality of
$m(n,k)$ since $v_2v_j>(v_2-1)(v_j+1)$.

\quad \mbox{\emph{Case 2: $v_j = n-i-1$}.} If $j=2$, as before:
$$1+(n-i-1)+(n-i+1)+\cdots+n\leq(n-i-1)+k-(n-i)-(n-i-1)-\cdots-(n-i_0)<k,$$
which contradicts (\ref{Ckn}). Thus $1<v_2<v_j$.

\quad \quad \emph{Subcase 2.1: $v_2 >2$}. Remplace $v_2$ by $v_2-1$
and $v_j$ by $v_j+1$ to get a contradiction since
$v_2v_j>(v_2-1)(v_j+1)$.

\quad \quad \mbox{\emph{Subcase 2.2: $v_2 =2$}.} Since
$k<\frac{n(n-1)}{2}$ (this inequality is crucial here), there exists
an integer $1<u<v_j$ such that $u$ is not one of the $v_i$. Consider
the greatest integer $l$ such that $v_l<u$. The minimality of
$m(n,k)$ is finally contradicted by replacing in the $r$-tuple
$\{2,v_l,v_j\}$ by $\{v_l+1,v_j+1\}$ since $2v_lv_j>(v_l+1)(v_j+1)$.
\end{preuve}

\begin{lem}\label{lem:ineg}Following the notations of lemma \ref{lem0}, we have:
\begin{equation}\label{squize}\frac{\Gamma(n+1)}{\Gamma(n-i_0)}\leq m(n,k) \leq \frac{\Gamma(n+1)}{\Gamma(n-i_0)} e^n,\end{equation}
where $\Gamma(n)$ is Euler's Gamma function.
\end{lem}
\begin{preuve}

By lemma \ref{lem0}, it is clear that
$\frac{\Gamma(n+1)}{\Gamma(n-i_0)}\leq m(n,k)$. For the second
inequality, let $j$ be the greatest integer such that $v_j<n-i_0$.
By definition of $i_0$:
$$k-1< n+(n-1)+\cdots+(n-i_0)+(n-i_0-1).$$
Since $(v_1,\ldots,v_r)$ satisfies (\ref{Ckn}):
$$1+\sum_{i=2}^j v_i+(n-i_0)+(n-i_0+1)+\cdots+(n-1)+n=k.$$
Hence:
$$\sum_{i=2}^j v_i<n-i_0-1\leq n$$
The arithmetic-geometric mean inequality and a brief study of the
function $x\mapsto \left(\frac{n}{x}\right)^x$ yield:
\begin{equation}\label{majoration}\prod_{i=2}^j v_i\leq \left(\frac{n}{j-1}\right)^{j-1}\leq e^n,\end{equation}
implying the result of the lemma.
\end{preuve}

\begin{prop}\label{prop2}For $n+1\leq k < \frac{n(n-1)}{2}$ we have:
\begin{equation}
\frac{1}{ne^n}\Gamma(n-i_0) \leq \C(n,k) \leq 2^n \Gamma(n-i_0)
\end{equation}
\end{prop}

\begin{preuve}
Combine proposition \ref{prop:bounds} and lemma \ref{lem:ineg}, and
use the fact that $\C(n,k)=\mathbb{P}(srec=k)/n!$.
\end{preuve}

\subsection{Case $\frac{n(n-1)}{2}\leq k \leq \frac{n(n+1)}{2}$}

Let $k$ be an integer such that $\frac{n(n-1)}{2}\leq k \leq
\frac{n(n+1)}{2}$ and $k \neq \frac{n(n+1)}{2}-1$. Let us prove that
the result of proposition \ref{prop2} holds in this case too (note
that we will not use lemma \ref{lem0}).

\begin{lem}\label{lembis:ineg}Let $(v_1,\ldots,v_r)$ be an $r$-tuple which realizes
the minimum $m(n,k)$ defined in (\ref{defmin}). Let $i_0=i_0(n,k)$
be the greatest integer such that:
$$k-1\geq n+(n-1)+\cdots+(n-i_0).$$ Then:
\begin{equation}\label{squize1}\frac{\Gamma(n+1)}{\Gamma(n-i_0)}\leq m(n,k) \leq \frac{\Gamma(n+1)}{\Gamma(n-i_0)} e^n,\end{equation}
where $\Gamma$ is Euler's Gamma function.
\end{lem}
\begin{preuve}It goes along the same lines as the proof of lemma \ref{lem0}.
If, for $0\leq i \leq i_0$, $n-i$ is equal to one of the $v_j$, then
we continue as in the previous subsection. If not, let $i$ be the
smallest integer such that $n-i$ does not appear in
$v_1,\ldots,v_r$. We can also assume that we are in the Subcase
$2.2$ of the proof of lemma \ref{lem0} (if not, we obtain a
contradiction as in the proof of lemma \ref{lem0}). Consequently,
$v_2=2$, there exists an integer $j$ such that $v_j<n-i$ and if
$1<u<v_j$ then $u$ is one of the $v_i$. In other words, all
positions but $n-i$ are records. For convenience, we introduce
$u=\frac{n(n+1)}{2}-k$, where $k=v_1+v_2+\cdots+v_r$, so that $3\leq
u \leq n$. Thus:
$$u=\frac{n(n+1)}{2}-k=n-i \quad \mbox{and} \quad m(n,k)=\frac{n!}{n-i}=\frac{n!}{u}.$$
By definition of $i_0$, Eq.(\ref{eqn:i0}), one gets:
\begin{equation}\label{expression}i_0(n,k)=\left[\frac{2n-1-\sqrt{4n^2+4n-8k+9}}{2}\right]=\left[n-\frac{1}{2}-\sqrt{2u+\frac{9}{4}}\right].
\end{equation}

\bigskip

Let us prove the first inequality in (\ref{squize1}). First note
that it is equivalent to $u\leq\Gamma(n-i_0)$. For $x>0$ define
$E(x)$ to be $x$ if $x \in \N$ and $[x]+1$ otherwise, so that
$n-[n-x]=E(x)$. In virtue of (\ref{expression}), the inequality
$u\leq\Gamma(n-i_0)$ is equivalent to:
$$u \leq\Gamma\left(E\left(\frac{1}{2}+\sqrt{2u+\frac{9}{4}}\right)\right).$$
This inequality is verified for $u=3$ and for all integers $u \geq
4$ one has:$$u\leq
\Gamma\left(\frac{1}{2}+\sqrt{2u+\frac{9}{4}}\right),$$ so that the
first inequality is proved.

\bigskip
Let us now prove the second inequality in (\ref{squize1}). For all
integer $u$ such that $3\leq u \leq n$, we have:
$$\Gamma(n-i_0)\leq\Gamma\left(\left[\frac{1}{2}+\sqrt{2u+\frac{9}{4}}\right]+1\right)\leq u
e^u\leq ue^n,$$ where the first inequality follows from the
properties of the $\Gamma$ function and the second one from the fact
that $u \leq n$.
\end{preuve}

Using proposition \ref{prop:bounds} and lemma \ref{lembis:ineg} we
finally extend proposition \ref{prop2} to:

\begin{prop}\label{prop3}For an integer $k$ such that $n+1\leq k \leq \frac{n(n+1)}{2}$ and $k \neq \frac{n(n+1)}{2}-1$ we have:
\begin{equation}
\frac{1}{ne^n}\Gamma(n-i_0) \leq \C(n,k) \leq 2^n \Gamma(n-i_0)
\end{equation}
\end{prop}

\subsection{Proof of the main theorem}
Let $n$ be an integer such that $n \geq 4$ and $x\in[0,1]$. Define
$k=k(n,x)=\left[x\frac{n(n+1)}{2}\right]$. Note that $3\leq k \leq
n$ if and only if $\frac{6}{n(n+1)}\leq x < \frac{2}{n}$. When
$x\geq \frac{2}{n}$, define $i_0(n,x):=i_0(n,k(n,x))$ as in
Eq.(\ref{eqn:i0}).

\begin{lem}\label{lemm}The following inequality holds:
$$\forall n \in \mathbb{N}, n\geq 4, \forall x \in \left[\frac{2}{n},1-\frac{2}{ n(n+1)}\right], \quad |n-i_0(n,x)-n\sqrt{1-x}|\leq 3.$$
\end{lem}

\begin{preuve}It can be deduced from Eq.(\ref{expression}) by using the fact that $k=\left[x\frac{n(n+1)}{2}\right]$.
\end{preuve}
\bigskip

\textbf{Proof of theorem \ref{thm2}}. Define $K_n(x)=\ln n
\left|\frac{\ln(\phi_n(x))}{n\ln(n)}-\sqrt{1-x}\right|$ (see
Eq.(\ref{defn:phin}) for the definition of $\phi_n$). Then:
$$\ee_{x \in [0,1]}K_n(x)\leq\ee_{x \in [0,\frac{6}{n(n+1)}]}K_n(x)+\ee_{x \in [\frac{6}{n(n+1)},\frac{2}{n})}K_n(x)+\ee_{x \in [\frac{2}{n},1]}K_n(x)$$
Denote the first term on the right-hand side of this inequality as
$A_n$, the second one as $B_n$ and the third one as $C_n$. We prove
that they are all bounded by a constant independent of $n$.
\begin{enumerate}
\item[(i)] We have:$$A_n=\ee_{x \in [0,\frac{6}{n(n+1)}]}\ln n
\left|\frac{\ln(n-1)!}{n \ln n}-\sqrt{1-x}\right|\leq
\widetilde{C}_1$$ by Stirling's formula.
\item[(ii)] By proposition \ref{prop1} we have:
\begin{align} B_n&\leq&&\ee_{x \in
[\frac{6}{n(n+1)},\frac{2}{n})}\left|\frac{\ln(n!)-\ln(\left[x\frac{n(n+1)}{2}\right]-1)-\ln
n}{n}-\sqrt{1-x}\ln n\right|    \notag\\
 & &&+\ee_{x \in
[\frac{6}{n(n+1)},\frac{2}{n})}\left|\frac{n\ln(2)+\ln(n!)-\ln(\left[x\frac{n(n+1)}{2}\right]-1)}{n}-\sqrt{1-x}\ln(n)\right|    \notag\\
&\leq&&\widetilde{C}_2+2\ee_{x \in
[\frac{6}{n(n+1)},\frac{2}{n})}\left|\frac{\ln(n!)}{n}-\sqrt{1-x}\ln(n)\right|  \notag\\
&\leq&& \widetilde{C}_3    \notag
\end{align}
\item[(iii)]By lemma \ref{lemm} for $\frac{2}{n}\leq x \leq
1-\left(\frac{5}{n}\right)^2$, one has $n\sqrt{1-x}\geq 5$ so that
$n-i_0(n,x)\geq n\sqrt{1-x}-3$ and $n\sqrt{1-x}-3\geq 2$. For
$1-\left(\frac{5}{n}\right)^2 \leq x < 1-\frac{2}{n(n+1)}$ one has
$n-i_0(n,x) \geq 2$. Note that for $x$ such that $x \geq
1-\frac{2}{n(n+1)}$ one has $\phi_n(x)=1$. Hence by proposition
\ref{prop3}:
\begin{align}
C_n&\leq&&\ee_{x \in [\frac{2}{n},1]}\left|\ln 2+\frac{\ln\left(\Gamma\left(n\sqrt{1-x}+3\right)\right)}{n}-\sqrt{1-x}\ln n\right|  \notag\\
& &&+\ee_{x \in [\frac{2}{n},1-\left(\frac{5}{n}\right)^2]}\left|-\frac{\ln n}{n}-1+\frac{\ln\left(\Gamma\left(n\sqrt{1-x}-3\right)\right)}{n}-\sqrt{1-x}\ln n\right|   \notag\\
& &&+\ee_{x \in [1-\left(\frac{5}{n}\right)^2,1-\frac{2}{n(n+1)})}\left|-\frac{\ln n}{n}-1-\sqrt{1-x}\ln n\right|   \notag\\
& &&+\ee_{x \in [1-\frac{2}{n(n+1)},1]}\left|\sqrt{1-x}\ln n\right| \notag\\
& \leq && \widetilde{C}_4,    \notag
\end{align}
and this concludes the proof.
\end{enumerate}
\fin

\section{Appendix}

In this appendix we show that our theorem \ref{thm1} is consistent
with Temme's result \cite{Temme} (see also \cite{Temme2}).

Let $m,n$ be positive integers such that $m \leq n$.
Define:\begin{equation}\label{eqn:defphi}\phi(u)=\ln((u+1)(u+2)\cdots(u+n))-m
\ln u.\end{equation} Let $u_1$ be the unique positive solution of
the equation $\phi'(u)=0$ (see \cite{Temme} for the proof that $u_1$
is unique). Let $t_1=m/(n-m)$ and $B=\phi(u_1)-n\ln(1+t_1)+m \ln
t_1$. Finally let
$g(t_1)=u_1^{-1}\left[m(n-m)/(n\phi''(u_1))\right]^{1/2}$.

\begin{thm}[Temme]\label{thm:temme}The relation
\begin{equation}c(n,m)\sim e^B g(t_1) \binom{n}{m}\end{equation}
holds uniformly for $1\leq m \leq n$ in the limit
$n\rightarrow\infty$.
\end{thm}

We show that our theorem \ref{thm1} is a consequence of theorem
\ref{thm:temme}. In other words, we deduce from Temme's formula the
fact that the coefficients $c(n,m)$ have a scaled asymptotic
behavior in the limit $n\rightarrow\infty$ with the ratio $m/n$
fixed, which is not clear a priori.

\begin{prop} Let $x$ be a real number such that $0<x<1$. Using the
previous definitions with $m=[nx]$, the following relation holds:
\begin{equation}\label{limi}\lim_{n\rightarrow\infty} \frac{\ln (e^B g(t_1) \binom{n}{[nx]})}{n \ln n}=1-x\end{equation}
\end{prop}

To prove this, the following lemma will be useful. It shows that
$u_1/n$ has a nice behavior for large $n$.

\begin{lem}\label{lem:encadr}For $n$ sufficiently large we have:
\begin{equation}\label{encad}\frac{x^2}{6(4/3-x)}\leq \frac{u_1}{n}\leq
\frac{x}{1-x}+\frac{1}{n}\end{equation}\end{lem}

\begin{preuve}Let
$f(u)=1/(u+1)+1/(u+2) +\cdots+1/(u+n)$. Recall that $u_1$ satisfies
the relation
\begin{equation}\phi'(u_1)=f(u_1)-\frac{[nx]}{u_1}=0.\end{equation} For the
second inequality in (\ref{encad}), note that:
$$\frac{n}{u_1+n}\leq \frac{1}{u_1+1}+\frac{1}{u_1+2}
+\cdots+\frac{1}{u_1+n} =\frac{[nx]}{u_1}\leq\frac{nx+1}{u_1},$$ so
that $u_1\leq (n(x+1/n))/(1-(x+1/n)) \leq nx/(1-x)+1$, where the
last inequality holds for $n$ sufficiently large.

For the first inequality in (\ref{encad}), let
$\alpha=x^2/(6(4/3-x))$ so that
$x/2+(1-x/2)\alpha/(\alpha+x/2)=3x/4$ and $p=[xn/2]$. Since the
function $u\mapsto u f(u)$ is increasing for positive $u$, it is
sufficient to show that $f(\alpha n) \leq (nx-1)/(\alpha n)$. Then:

\begin{align}f(\alpha n)&=&&\frac{1}{\alpha n+1}+\frac{1}{\alpha n+2}
+\cdots+\frac{1}{\alpha n+n}\leq \frac{p}{\alpha
n}+\frac{n-p}{\alpha n+p}\leq\frac{xn/2}{\alpha
n}+\frac{n-(xn/2-1)}{\alpha n + xn/2-1}\nonumber\\\nonumber
&\leq&&\frac{x}{2\alpha}+\frac{1-x/2-1/n}{\alpha+x/2-1/n} = \frac{3x}{4\alpha}-\frac{1-x/2}{\alpha+x/2}+\frac{1-x/2-1/n}{\alpha+x/2-1/n}\\
&\leq&&\frac{x}{\alpha}-\frac{1}{n \alpha},\end{align} where the
last inequality holds for $n$ sufficiently large.
\end{preuve}

\textbf{Proof of proposition 5.2.} Let us show that only the term
$\phi(u_1)$ entering the expression for $B$ provides the dominant
contribution on the right-hand side of Eq.(\ref{limi}). More
precisely:
\begin{enumerate}
\item[-]Using Stirling's formula, it is easy to see that $\lim_{n \rightarrow
\infty}\frac{\ln \binom{n}{[nx]}}{n \ln n}=0$.
\item[-]To show that $\lim_{n \rightarrow
\infty}\frac{\ln g(t_1)}{n \ln n}=0$, it is sufficient to show that
$\lim_{n \rightarrow \infty} \frac{\ln \phi''(u_1)}{n \ln n}=0$. It
is convenient to write $\phi''(u_1)$ as
$\phi''(u_1)=\psi'(u_1+n+1)-\psi'(u_1+1)+[nx]/u_1^2$, where $\psi$
is the psi function. Lemma \ref{lem:encadr} and the fact that
$\psi'(x)\sim 1/x$ for large $x$ (see e.g. \cite{Abramowitz}) give
the result.
\item[-]It is clear that $\lim_{n \rightarrow
\infty} \frac{\ln |-n\ln(1+t_1)+[nx] \ln t_1|}{n \ln n}=0$.
\item[-]Finally, observe that $\phi(u_1)=\ln \Gamma(u_1+n+1)-\ln
\Gamma(u_1+1)-[nx]\ln u_1$. Stirling's formula and the fact that
$\ln(u_1+n+1)=\ln(u_1+1)+o(\ln n)$ (which is a consequence of lemma
\ref{lem:encadr}) show that $\lim_{n \rightarrow
\infty}\frac{\phi(u_1)}{n \ln n}=1-x$. This concludes the proof.
\end{enumerate}
\fin

Thus we have reproduced the scaled asymptotic behavior of
$c(n,[nx])$ using Temme's result. However, it seems difficult to
also reproduce by this means the error estimate stated in theorem
\ref{thm1}. To this end, it would be necessary to give a more
precise asymptotic behavior of $u_1$.

\section*{Conclusion}
We have studied the asymptotic behavior of the integers $c(n,k)$
(respectively $\C(n,k)$) equal to the number of elements of $\S_n$
having $k$ records (respectively for which the sum of the positions
of their records are $k$) by using a probabilistic argument. One can
note that these integers can be defined outside of any combinatorial
background since $c(n,k)$ appears as the coefficient of $q^k$ in the
polynomial $q(q+1)\cdots(q+n-1)$ and $\C(n,k)$ appears as the
coefficient of $q^k$ in the polynomial
$q(q^2+1)(q^3+2)\cdots(q^n+n-1)$. Thus studying the asymptotic
behavior of these numbers seems delicate, but the probabilistic
interpretation gave us a convenient formula defining these integers.
Surprisingly, the scaled asymptotic behavior of these rather
complicated numbers can be described by a remarkably simple
function.

\section*{Acknowledgments}

I am deeply indebted to Richard Stanley for his helpful advice as
well as to Philippe Biane for useful discussions and for sharing the
idea which lead to the proof of the main theorem. I thank Herbert S.
Wilf for useful comments and for pointing out additional references.

\section*{Note added}

Recently, based on the result obtained in the present paper, the
statistic `sum of the position of records' has also been considered
in the case of the geometric law in \cite{Geo3}.


\begin{thebibliography}{99}

\bibitem{Renyi}
A. Rényi, Théorie des éléments saillants d'une suite d'observations
, \emph{Ann. Fac. Sci. Univ. Clermont-Ferrand} No. \textbf{8} 1962
7-13.

\bibitem{Jordan}
C. Jordan, The calculus of finite differences, 2nd ed., Chelsea,
1947.

\bibitem{MoserWyman}
L. Moser, M. Wyman, Asymptotic development of the Stirling numbers
of the first kind, \emph{J. London Math. Soc.} \textbf{33}: 133-146,
1958.

\bibitem{Wilf}
H. S. Wilf, The asymptotic behavior of the Stirling numbers of the
first kind, \emph{J. Combin. Theory Ser. A}, \textbf{64}:
344-349,1993

\bibitem{Temme}
N. M. Temme, Asymptotic estimates of Stirling numbers, \emph{Studies
in Applied Math}, \textbf{89}: 233-243, 1993.

\bibitem{Hwang}
H. Hwang, Asymptotic expansions for the Stirling numbers of the
first kind, \emph{J. Combin. Theory Ser. A}, \textbf{71(2)}:
343-351, 1995.

\bibitem{Temme2}
R. Chelluri, L.B. Richmond, N.M Temme, Asymptotic estimates for
generalized Stirling numbers, \emph{Analysis}, \textbf{20}: 1-13,
2000.

\bibitem{Wilf1995}
H. S. Wilf, On the outstanding elements of permutations,
\verb"http://www.math.upenn.edu/~wilf/website/outstelmts.pdf", 1995.

\bibitem{Geo0}
H. Prodinger, Combinatorics of geometrically distributed random
variables: Left-to-right maxima. \emph{Discrete Mathematics},
\textbf{153}: 253-270, 1996.


\bibitem{Geo}
A. Knopfmacher, H. Prodinger, Combinatorics of geometrically
distributed random variables: value and position of the rth
left-to-right maximum, \emph{Discrete Mathematics}, \textbf{226}:
255-267, 2001.

\bibitem{Geo2}
H. Prodinger, Combinatorics of geometrically distributed random
variables: Value and position of large left-to-right maxima,
\emph{Discrete Mathematics} \textbf{254}: 459-471, 2002.

\bibitem{Geo3}
H. Prodinger, Records in geometrically distributed words: sum of
positions, \emph{Appl. Anal. Discrete Math.} \textbf{2}: 234-240,
2008.

\bibitem{MW}
A. N. Myers, H. S. Wilf, Left-to-right maxima in words and multiset
permutations, \emph{Israel J. Math.} \textbf{166}: 167-183, 2008
[arXiv:math/0701078]

\bibitem{Glick}
N. Glick, Breaking records and breaking boards. \emph{American
Mathematical Monthly}, \emph{85}:2-26, 1978.

\bibitem{PR}
R. P. Stanley, \emph{Private communication}, Unpublished.

\bibitem{ENUM}
R. P. Stanley, \emph{Enumerative Combinatorics, vol. 1}, Cambridge,
England: Cambridge University Press (1999).

\bibitem{Abramowitz}
M. Abramowitz, I. A. Stegun, \emph{Handbook of Mathematical
Functions with Formulas, Graphs and Mathematics Tables}, National
Bureau of Standards, Applied Series 55, U.S. Government Printing
Office, Washington, D.C. , 1964.



\end{thebibliography}
\end{document}